\documentclass{amsart}

\newtheorem{thm}{Theorem}[section]

\newtheorem{lem}[thm]{Lemma}

\theoremstyle{definition}

\theoremstyle{remark}

\numberwithin{equation}{section}

\newcommand{\R}{\mathbb R}

\newcommand{\si}{\sigma}
\newcommand{\rt}{\rightarrow}

\newcommand{\w}{\mathcal{W}}
\newcommand{\g}{\hat{g}}
\newcommand{\la}{\lambda}

\begin{document}

\title[ Weyl curvature and the Euler characteristic]{Weyl curvature and the Euler characteristic in dimension four }
\author{harish seshadri}
\address{department of mathematics,
Indian Institute of Science, Bangalore 560012, India}
\email{harish@math.iisc.ernet.in}



\begin{abstract}
We give lower bounds, in terms of the Euler characteristic, for
the $L^2$-norm of the Weyl curvature of closed Riemannian
4-manifolds. The same bounds were obtained by Gursky, in the case
of positive scalar curvature metrics.
\end{abstract}
\maketitle
\section{Introduction}
Let $M$ be a smooth closed oriented $4$-manifold and let $C=[g]:=
\{ fg : f \in C^\infty (M)$ and $f>0 \}$ be a conformal class of
metrics on $M$. An important numerical invariant associated to
$C$ is the {\it Weyl constant} $\w (C)$. The Weyl constant is
defined by
$$ \w (C)=\int_M \vert W_g \vert_g ^2dV_g,$$
where $g$ is {\it any } metric in $C$ and $W$ is the Weyl tensor
of $g$. Since the vanishing of the Weyl tensor is equivalent to
the conformal flatness of $g$, one can regard $\w (C)$ as a
quantitative measure of the lack of conformal flatness.

As the existence of a conformal class with prescribed value of
$\w$ is a diffeomorphism invariant, one can try to relate $\w$ to
standard topological invariants. In fact, in dimension 4 one has

\begin{thm} \label{one}{\rm ({\bf Gursky} ~\cite{gur1})}
Let $(M,g)$ be a closed oriented Riemannian 4-manifold. If $g$ has
positive scalar curvature, then $$\int_M \vert W \vert ^2 \ge 8
\pi^2 (\chi (M) -2).$$ Equality holds if and only if $g$ is
conformal to an Einstein metric $h$ with $s_h \ Vol_h^ { \ {\frac
{1}{2}}} =8 \pi \sqrt 6$, where ``$s$" denotes scalar curvature.
\end{thm}
Note that $8 \pi \sqrt 6$ is the {\it Yamabe constant} of the
standard metric on $S^4$. Hence the results of Schoen ~\cite{sch}
imply that $(M,g)$ is conformally equivalent to $S^4$ in the case
of equality above.

As a corollary of Theorem \ref{one}, one obtains
\begin{thm}\label{two}{\rm ({\bf Gursky} ~\cite{gur1})}
Let $(M,g)$ be a closed oriented Riemannian 4-manifold. If $g$ is
conformally flat and has positive scalar curvature, then $\chi
(M) \le 0$ unless $(M,g)$ is conformally equivalent to the round
4-sphere.
\end{thm}

Theorems \ref{one} and \ref{two} were proved by Gursky in
~\cite{gur1}(Theorem \ref{one} is not stated as such but is
contained in the proofs). In the first part of this paper
(Section 2) we give a simple, geometric proof of these results
using ``stereographic" projection. As noted by Gursky, the proofs
of these results would be relatively straightforward if one were
to assume the existence of a Yamabe metric in every conformal
class. However, the known proof of existence of a Yamabe metric in
dimension 4 uses the hard and deep Positive Mass Theorem of
Schoen and Yau. Hence, in order to make the proofs ``elementary",
we try to avoid the use of a Yamabe metric and use it only for
the case of equality in Theorem \ref{one}.

In the second part (Section 3) of the paper, we prove a version of
Theorem \ref{one} for {\it nonpositive} scalar curvature metrics.

\begin{thm}\label{thr}
Let $(M,g)$ be a closed oriented Riemannian 4-manifold. If $ s +
c \vert W \vert  \ge 0 $ for some $c >0$ and there is a metric $h$
conformal to $g$ with $ \int_M s_h \le 0$, then $$ \ \int_M \vert
W \vert ^2 \ge \frac {8 \pi^2}{1+ \frac {c^2}{24}} \chi (M).$$
Equality holds if and only if $g$ is an Einstein metric with $s +
c \vert W \vert \equiv 0$.
\end{thm}

Let us note that the hypotheses (and the conclusion) in the above
theorem are dependent only on the conformal class of the metric
$g$.

It should be mentioned that different (and far subtler) sharp
lower bounds for $\int_M \vert W \vert ^2$ were obtained by
Gursky in ~\cite{gur2} (for positive scalar curvature metrics,
under the assumption of non-zero first or second Betti number)
and in ~\cite{gur3} (for negative scalar curvature metrics, under
the assumption of the existence of a conformal vector field).

Our strategy for proving Theorems \ref{one} and \ref{two} is to
use the ``stereographic projection" of $(M,C)$. This gives us a
complete noncompact asymptotically flat scalar-flat 4-manifold
$(\hat M, \hat g)$. The two main points for us are: First, under
this passage, the Weyl invariant does not change. Second, the
scalar-flatness and asymptotic flatness simplify the
Chern-Gauss-Bonnet formula for balls in $(\hat M, \hat g)$
considerably. Unfortunately, it is is not clear how to extend
this method to dimensions beyond 4 since we crucially use the
specific form that Chern-Gauss-Bonnet takes in this dimension.

In Section 3 we prove Theorem \ref{thr}. We use Yamabe metrics in
this case. It should be possible, with some extra effort, to give
a proof using stereographic projections but we do not do pursure
this approach here. \\
{\bf Acknowledgements} \\
It's a pleasure to thank Mike Anderson for helpful discussions.

\section{Stereographic projection and the Weyl constant}
For rest of this section we assume that $(M,g)$ is a closed
oriented Riemannian 4-manifold with positive scalar curvature.
Fix $p \in M$ and let $G$ denotes the Green's function of the
conformal Laplacian $L= 6\triangle - s$ at $p$. Since $s >0$, $G$
exists and is positive. Also $\g = G^2g$ is a complete,
scalar-flat, asymptotically flat metric on $\hat M :=M- \{ p \}$
(cf. ~\cite{lp}). $(\hat M, \hat g)$ is sometimes referred to as
the "stereographic projection" of $(M, [g])$. Let $S_r$ and $B_r$
denote the sphere and closed ball of radius $r$ at $p$ in $(M,g)$.

Lemmas (\ref{pc}) and (\ref{bd}) will imply that the boundary
integral in Chern-Gauss-Bonnet applied to certain large domains in
$(\hat M, \hat g)$ will give the same value as for balls in flat
$\R^4$. The domains we consider are the complements of $B_r$ in
$M$. Forms of the Gauss-Bonnet theorem for asymptotically flat
manifolds have been described in ~\cite{gil} and ~\cite{and0}.
For the specific result that we need and for the sake of
completeness, we give the computation of the boundary integrals
in detail.

In the next lemma the principal curvatures are with respect to the
{\it inward} pointing normal of $S_r \subset \hat M$.

\begin{lem}\label{pc}
If $\hat \la_r$ is a principal curvature of $S_r$ with respect to
$\g$, then $\hat \la(r) = -r + O(r^2)$ as $r \rt 0$.
\end{lem}
\begin{proof}
In what follows, hats will denote quantities defined with respect
to $\g$. The second fundamental form $\hat B$ of $S_r$ is related
to $B$ by
$$ \hat B = G \ B +  \frac {\partial G}{\partial r} g,$$
where we have used standard formulas for conformal changes. Hence
we have the following equations for the shape operator $S$, which
is given by $B(X,Y)=g(S(X),Y)$, and the principal curvatures,
which are the eigenvalues of $S$: $$\hat S = G^{-1} S +
G^{-2}\frac {\partial G}{\partial r} I, \ \ \ \hat \la_r = G^{-1}
\la_r + G^{-2}\frac {\partial G}{\partial r}.$$ Now let $\{ x^i
\}$ denote {\it conformal normal coordinates } at $p$, as defined
in ~\cite{lp}. If $r=d(x,p)$, then we have
\begin{equation}\label{gre}
 G(x) = r^{-2} + A+O''(r) \ \ {\rm as} \ r \rt 0,
\end{equation}
where $f=O''(r^k)$ means $f=O(r^k)$, $\nabla f =O(r^{k-1})$ and
$\nabla ^2 (f) =O(r^{k-2})$. We do not use this information but we
note that $A$ is a multiple of the {\it mass} \ of the
asymptotically flat manifold $(\hat M, \ \g)$. From the above
expression we get $\frac {\partial G}{\partial r} = -2r^{-3}+
O(1)$. Finally
\begin{equation}\label{mai}\notag
\hat \la_r  \ = \ G^{-1} \la_r + G^{-2}\frac {\partial
G}{\partial r}
        \ = \ \frac{ r^{-1} + O(r)}{r^{-2} + A + O(r)} + \frac{-2r^{-3} + O(1)}{r^{-4}
        +O(r^{-2})}
       \ = \ -r + O(r^2) \notag\\
\end{equation}
In the second equality we have used the well-known (see
~\cite{gl}, for instance) and easily verified fact that $\la_r
=r^{-1} +O(r)$ on any Riemannian manifold .
\end{proof}
\medskip
The Chern-Gauss-Bonnet formula for a manifold with boundary $N$
states  (see ~\cite{che} and also ~\cite{and}) that
\begin{equation}\label{che}
8 \pi^2  \chi (N) = \int_N (\vert W \vert ^2 -\frac {1}{2} \vert z
\vert ^2  \ + \frac {1}{24} s^2) \ - 4 \int_{\partial N} \prod_1^3
\la_i \ - \ \int_{\partial N} \Sigma_{\sigma \in S_3} K_{\sigma_1
\sigma _2} \la_{\sigma_3}
\end{equation}
Here $W$, $z=ric -\frac {s}{n}g$ and $s$ are the Weyl, trace-free
Ricci and scalar curvature, respectively, $K$ denotes sectional
curvature and $\la_i$ the principal curvatures of $\partial N$.
Let us denote by $I^1_r= \int_{S_r} \prod_1^3 \la_i $ and $I^2_r =
\int_{S_r} \Sigma_{\sigma \in S_3} K_{\sigma_1 \sigma _2}
\la_{\sigma_3}$ the two boundary integrals in the above formula
applied to $(M_r, G^2 \ g):= (\ M - Int \ B_r, G^2 \ g))$.
\begin{lem}\label{bd}
$lim_{r \rightarrow 0} \ I^1_r = -2 \pi^2$ and $lim_{r \rightarrow
0} \ I^2_r=0$.
\end{lem}
\begin{proof}
If $dA_r$ denotes the volume form of $S_r$ in $(M_r,g)$, then $
\hat {dA}_r = G^3dA_r= (r^{-6}+O(r^{-4}))dA_r$. Now
$$ I^1_r= \int_{S_r} \prod_1^3  \hat \la_i \hat {dA}_r
=\int_{S_r} (-r^3 + O(r^4))(r^{-6} +O(r^{-4})dA_r,$$ where we have
used (\ref{mai}) in the last equation. Since $Vol (S_r) =O(r^3)$,
$$lim_{r \rightarrow 0} \ I^1_r = lim _{r \rightarrow 0} \ -r^{-3}Vol
(S_r)= -2 \pi^2.$$ The last equation above can be easily seen by
using normal coordinates. As for $I^2_r$, it is clear from
(\ref{gre}) that for $r$ small enough, $\vert K \vert \le 1$ on
$S_r$. Hence the integrand (with respect to $dA_r$) in $I^2_r$ is
of $O(r^{-2})$ and $I^2_r \rightarrow 0$, as above.
\end{proof}
Now we come to the proof of Theorem \ref{one}.
\begin{proof}
Applying the Chern-Gauss-Bonnet theorem to $(M_r, G^2g)$, setting
$\hat s=0$ and getting rid of the $\vert \hat z \vert ^2$ term,
we get
$$8 \pi^2 \chi (M_r) \le \int_{M_r} \vert \hat W \vert ^2d \hat V + I^1_r +
I^2_r .$$ From the conformal invariance of $W$, we have
$\int_{M_r} \vert \hat W \vert ^2 d \hat V=  \int_{M_r} \vert W
\vert ^2 dV \rt \w (M)$ as $r \rt 0$. By Lemma \ref{bd} and
(\ref{che}) we have
$$8 \pi^2 \chi (M - \{ p \})= 8 \pi ^2 \chi(M_r) \le \w (M) + 8 \pi^2.$$
Since $ \chi (M) = \chi (M - \{ p \}) + 1$, we finally get $ \w
(M) \ge 8 \pi^2 (\chi (M)-2).$

Now suppose that
\begin{equation}\label{eq}
 \w (M) =8 \pi^2 (\chi (M)-2).
\end{equation}
Let $C$ denote the conformal class of $g$. Let $h \in C$ be a
Yamabe metric, i.e, a metric minimizing the total scalar
curvature functional $E$
$$ \tilde g \rt \ E(\tilde g) = \frac {\int_M s_{\tilde g} dv_{\tilde g}}{Vol
(\tilde g)^ {\frac {1}{2}}}, \ \ \ \tilde g  \in C.$$ The
existence of $h$ is guaranteed by ~\cite{sch}. $h$ has constant
scalar curvature, which implies that
\begin{equation} \label{con}
\int_M s_h^2 = \frac {( \int_M s_h)^2}{Vol (h)}
\end{equation}
Moreover, by Aubin, the infimum of $E$ cannot be greater than the
value of $E$ on the round sphere:
\begin{equation}\label{aub}
\frac { \int_M s_h}{Vol (h)^{\frac {1}{2}}} \le 8 \pi \sqrt 6.
\end{equation}
Combining (\ref{eq}) and the Chern-Gauss-Bonnet formula for $M$,
we get
$$-\frac {1}{2} \vert z_h \vert^2+ \frac {1}{24} \int_M s_h^2 -16 \pi^2 =0,$$
By (\ref{con}) and (\ref{aub}) we see that the sum of the last two
terms above is nonpositive. Hence we must have $z_h=0$, i.e., $h$
is Einstein, and also $\frac {1}{24} \int_M s_h^2 =16 \pi^2$.
\end{proof}

Now we come to the proof of Theorem \ref{two}.
\begin{proof}
We assume that $(M,g)$ is conformally flat, i.e. $W=0$. If $g$ has
positive scalar curvature and $\chi (M)
>0$, the we claim that $\chi (M)=2$. This is because $\chi (M)=
2-2 \beta_1 + \beta_2$, by Poincare Duality. However $\beta_2=0$
by the Bochner formula for harmonic 2-forms on $(M,g)$.

Since $\chi (M)=2$ and $g$ is conformally flat, we can appeal to
the theorem above and conclude that a Yamabe metric $h$ in $[g]$
is Einstein. Since $W_h=0$ it would follow that $h$ is of constant
(positive) sectional curvature and by orientability, it would
follow that $(M,h)$ is isometric to $(S^4,g_0)$ and we would be
done. However, we give a different proof which avoids the
existence of a Yamabe metric: First, $\chi (M_r)=\chi (M - \{ p
\})=1$. Again by applying Lemma \ref{bd} and (\ref{che}) to $M_r$
and letting $r \rt 0$, we get $\int_{\hat M}\vert \hat z \vert ^2
=0$. Combining this with $\hat s=0$ and $\hat W=0$, we see that
$( \hat M, \hat g)$ is a complete noncompact flat 4-manifold.

The Bieberbach theorem combined with the fact that $\chi ( \hat M)
= \chi (M_r)\neq 0$ imply that $\hat M$ is simply-connected and
hence isometric to flat $\R ^4$. It then follows that $(M,g)$ is
conformally equivalent to $S^n$: If $G_0$ is the Green's function
(with singularity at the north pole) for the conformal Laplacian
on $(S^4, g_0)$, then $ G^{-2}_0G^2g$ is a metric of constant
curvature $1$ on $M - \{ p \}$ (we have identified $\hat M$ with
$ \R ^4$), which extends to a smooth metric conformal to $g$ on
$M$.
\end{proof}

\section{Nonpositive scalar curvature and the Weyl constant}

Here we prove Theorem \ref{thr}. So assume that $s +c \vert W
\vert \ge 0$ for $g$.

\begin{lem} \label{pos}
For any metric $h$ in $[g]$ we have $\int_M (s_h+c \vert W_h
\vert ) \ge 0$.
\end{lem}
\begin{proof}
Let us introduce, following Gursky and LeBrun ~\cite{gl2}, the
{\it modified scalar curvature} $\sigma_g = s_g + c\vert W_g
\vert _g$. Under a conformal change $g \rightarrow \tilde g=u^2
g$, the modified scalar curvature transforms (with our convention,
the Laplacian $\triangle = \frac{d^2}{dx^2}$ on $\R$) by
\begin{equation}
\sigma _g \rightarrow \sigma _{\tilde g}  = u^{-2} \sigma _g -6
u^{- 3} \triangle u. \label{tra}
\end{equation}

We also have the functional
\begin{equation}
E_{\si} (g) = \int_M \sigma _g dV_g.
\end{equation}

When we restrict $E_\si$ to a conformal class we get an operator
$L$ on $C^\infty(M)$ defined by $L(u)= E_\si(u^2g)$ or
$$L(u)=-6 \triangle u + \sigma _g u.$$
Let $<,>$ denotes the $L^2$ inner product on $C^\infty(M)$. and
let
$$\lambda= inf _{{f \in W^{1,2}} \atop {\Vert f \Vert _2=1}}<Lf,f>$$
and $u$ be the corresponding eigenfunction. We note that since
$\sigma _g$ is, in general, Lipschitz continuous but not smooth
(at the zero locus of $\vert W \vert$), the best regularity we
can obtain for $u$ is that $u \in C^{2, \alpha}$ for any $0<
\alpha <1$. This is sufficient for our purposes. By the minimum
principle $u>0$ and by definition, $u$ satisfies
\begin{equation}
 L(u)= \lambda u. \label{pos}
\end{equation}
{\it Claim}: If $g_0=u^2g$, then  $\sigma_{g_0} \ge 0$. \\
{\it Proof}: By (\ref{tra}) and (\ref{pos}), we see that
$$ \sigma _{g_0}= u^{- 3} L(u) = \lambda u^ {- 2}.$$
Hence $\sigma_{g_0}$ has a fixed sign. Suppose $\sigma_{g_0} < 0$.
Then (\ref{tra}) would imply that $\triangle u \le 0$. Hence, by
the minimum principle, $u$ would be constant. But this would
contradict (\ref{tra}). This proves the Claim.

Suppose that $h = f^2 g$. Let $ h_0= \Vert f \Vert_2 ^{-2}f^2 g$.
Then $ E_{\si} (h_0) \ge E_{\si} (g_0) >0$. Since $ E_{\si}(h) =
\Vert f \Vert_2 ^2 E_{\si}(h_0)$, the lemma is proved.
\end{proof}

Now let $h$ be a metric of constant scalar curvature in $[g]$.
This exists by the solution to the Yamabe problem. Note that since
we have assumed that $\int_M s_g dV_g \le 0$, the Yamabe metric
has nonpositive scalar curvature. We work with $h$ for rest of the
proof. By the Lemma \ref{pos}, we have $\int_M(s + c \vert W
\vert) \ge 0$. Hence
\begin{align}
 \int_M s^2 = Vol^{-1} (\int_M s)^2  & \le c^2 Vol^{-1} (\int_M \vert W \vert
 )^2 \notag \\
 & \le c^2 \int_M \vert W \vert ^2 \notag ,
\end{align}
where we have used the Cauchy-Schwartz inequality at the last
step.

 Combining this with the Chern-Gauss-Bonnet formula, we are done.
If equality holds in Theorem \ref{thr}, we must have $s_h +c \vert
W_h \vert \equiv 0$ and $h$ must be Einstein. Again referring to
(\ref{tra}), we see that $h$ must be a constant multiple of $g$.
Hence $\sigma_g \equiv 0$ and $g$ must be Einstein. \hfill
$\square$
\medskip

For $M$ which do not admit positive scalar curvature metrics, it
would be interesting to estimate (in terms of the topology of
$M$) the smallest $c$ such that $\sigma \ge 0$ for some $g$.


\end{document}